\hsize=15cm
\vsize=22cm
\hoffset 1cm

\def\t{\textstyle}
\vglue 2cm
\centerline{\bf SELF-DUALITY IN DIMENSIONS $2n>4$: EQUIVALENCE OF VARIOUS}
\centerline{\bf DEFINITIONS AND AN UPPER BOUND FOR ${\bf p}_2$}
 \vskip 1cm
\centerline{Ay\c se H\"umeyra Bilge }
\vskip .5cm
\centerline{TUBITAK-Marmara Research Center}
\centerline{Research Institute for Basic Sciences}
\centerline{Department of Mathematics}
\centerline{P.O.Box 21, 41470 Gebze-Kocaeli, Turkey}
\centerline{e-mail:bilge@yunus.mam.tubitak.gov.tr}
\vfill
\centerline{\bf Abstract}
\baselineskip 16pt
\vskip .5cm

We show that the self-duality defined in [Trautman, Int.J.Theor.Phys.,{\bf
16},561 (1977)] is equivalent
to strong self-duality defined in [Bilge, Dereli and Kocak,
Lett.Math.Phys., {\bf 36}, 301-309, (1996)] and we obtain
an upper bound on $\int (p_2)^{n}$, where $p_2$ is the second Pontrjagin
class of an $SO(N)$ bundle over a $8n$ dimensional manifold.

 \vskip 2cm
\eject
\vskip 1cm
\baselineskip 22pt
\noindent
{\bf 1. Introduction.}

The self-duality of a 2-form in four dimensions is defined to be the Hodge
duality [1].
The search for a
meaningful definition of self-duality in higher dimensions as a
map on 2-forms, has been studied by various authors [2-7].
A definition of self-duality based on an eigenvalue criterion has been
recently proposed
in [8] where an overview of previous results is also given.
According to this definition, the
``strongly self-dual" 2-forms in $2n$ dimensions form  an $n^2-n+1$ dimensional
submanifold ${\cal S}_{2n}$ of the skew-symmetric $2n\times 2n$ matrices,
determined by the condition that the absolute values of the
eigenvalues coincide.  It was also shown that this definition is equivalent
to the self-duality condition used by Grossman et.al. [6], namely
$F\wedge F$ be self-dual in the Hodge sense.
Here, in Section 2 we show
that the self-duality condition given by Trautman [3]
$\omega^{n-1}=k*\omega$, where $*$ denotes the Hodge dual,  is also
equivalent to ``strong self-duality".

In [8], it has also  been shown that for an $SO(N)$ bundle,
the strongly self-dual forms saturate an upper bound for the integral of
the first Pontrjagin class  $p_1$. Here we obtain an upper bound for
the integral of the second Pontrjagin class $p_2$ for an $SO(N)$ bundle, which
is again saturated for strongly self-dual forms. This upper bound is however
gauge independent only under an ``orthonormality" assumption (Eq.3.4).

\vskip 1cm

\noindent
{\bf 2. Equivalence of various definitions of self-duality in $2n$ dimensions.}
\vskip .2cm
\noindent
{\bf 2.1. Preliminaries.}

Let $M$ be a $2n$ dimensional manifold, $E$ be a vector bundle over
$M$ with standard fiber $R^N$ and structure group $G$, and let {\bf g} be a
linear representation of  the Lie algebra of G. If $A$ is the {\bf g} valued
local connection 1-form, the local curvature 2-form is defined as
$F=dA-A\wedge A$ and the Bianchi identites are
$dF+F\wedge A-A\wedge F=-0.$
 The local curvature
2-form depends on the trivilization of the bundle, but its invariant
polynomials $\sigma_i$ defined by
${\rm det}(I+tF)=\sum_{i=0}^n\sigma_it^i$
are invariant of the local trivialization. The $\sigma_i$ 's
are closed $2i$-forms,
defining
deRham cohomology classes in $H^{2i}$, proportional to the Chern classes
${\bf c}_i$ (or the
Pontrjagin classes ${\bf p}_{i/2}$) of the bundle $E$. In addition for an
$SO(N)$ bundle, the
square root of the determinant of $F$ (which is a ring element) defines the
Euler class $\chi$ [9]. The $\sigma_i$'s can also be written as
linear combinations of $trF^i$ [10], where $F^i$
means the product of the
matrix $F$ with itself $i$ times, with the wedge multiplication of the entries.
In order to avoid proportionality constants, we will
work with the quantities $\sigma_i$'s instead of the Chern or Pontrjagin
classes.
In the following we consider real $SO(N)$ bundles, and we recall that
all principal minors of $F$  are
skew-symmetric matrices,
hence their determinants are perfect squares, and  the $\sigma_{2i}$'s are sums
of squares of $i$-fold
products of the entries of $F$.

We shall use the notation, $(a,b)$ where $a,b$ are forms on $M$
to represent the inner product
in on the exterior algebra $\Lambda M_x$, $x\in M$, while $\langle A,B\rangle$
where $A, B$ are Lie algebra
valued forms will denote the inner product in $\Lambda M_x\otimes E_x$, $x\in
M$,
hence both expressions are functions defined on the manifold.

\vskip .5cm
\noindent
{\bf 2.2. Equivalence of strong self-duality and Trautman's definition.}

Let $\omega_{ij}$ be the components of a 2-form in $2n$ dimensions with respect
to some local orthonormal basis. In the following, we shall denote the 2-form
$\omega$ and the skew-symmetric matrix consisting of its components with
respect to some orthonormal basis by the same symbol.  Let
$\pm i\lambda_{k}$, $k=1\dots n$ be the eigenvalues of the skew-symmetric
matrix $\omega$. The invariant
polynomials $s_{2i}$, of $\omega$ can be expressed in terms of the elementary
symmetric functions of the $\lambda_i^2$'s.
The inner products
$(\omega^i,\omega^i)$ and the $s_{2i}$'s are related as follows.
$$\eqalignno{
(\omega,\omega)=&s_2=\lambda_1^2+\lambda_2^2+\dots
+\lambda_n^2,\cr
\t{1\over (2!)^2}(\omega^2,\omega^2)=&
s_4=\lambda_1^2\lambda_2^2+\lambda_1^2\lambda_3^2+\dots
                +\lambda_{n-1}^2\lambda_n^2 \cr
\t{1\over (3!)^2}
(\omega^3,\omega^3)=&s_6=\lambda_1^2\lambda_2^2\lambda_3^2
                +\lambda_1^2\lambda_2^2 \lambda_4^2+\dots
                +\lambda_{n-2}^2\lambda_{n-1}^2 \lambda_n^2\cr
\dots&\quad \dots\cr
\t{1\over (n!)^2}(\omega^n,\omega^n)=
{1\over (n!)^2}\mid*\omega^{n}\mid^2=&s_{2n}=
               \lambda_1^2\lambda_2^2\dots       \lambda_n^2&(2.1) \cr}$$

  We define the weighted elementary
symmetric polynomials by $q_i=s_{2i}/{n\choose i}$. We have the inequalities
[11]
$$q_1\ge q_2^{1/2}\ge q_3^{1/3}\ge \dots \ge q_n^{1/n},\quad\quad
q_{r-1}q_{r+1}\le q_r^2, \quad 1\le r<n\eqno(2.2)$$
and the equalities hold iff all the $\lambda_i$'s are equal.
In [8] we have adopted this case as a definition.
\vskip .3cm

\noindent
{\bf Definition 2.1 }
Let $\omega$ be a  2-form
in $2n$ dimensions,  $\pm i\lambda_k$, $k=1,\dots,n$ be its
eigenvalues and $\eta$ be the the square root of the determinant of $\omega$,
with a fixed choice of sign.
Then $\omega$ is called {\it strongly self-dual} ({\it strongly anti
self-dual}) if $\mid\lambda_1\mid=\mid\lambda_2\mid=\dots=\mid\lambda_n\mid$,
and $\eta>0$ ($\eta<0$).
\vskip .2cm

The strong self-duality condition is
equivalent to
the matrix equation $\omega^2+\lambda^2 I=0$, where $I$ is the identity
matrix, and $\lambda^2={1\over 2n}{\rm tr}\omega^2$. This definition gives
quadratic equations for the $\omega_{ij}$'s, hence the strong
self-duality condition determines a nonlinear set. In four dimensions, the
strong self-duality coincides with usual Hodge duality, more precisely, the
matrices satisfying $\omega^2+\lambda I=0$
consist of the union of the usual  self-dual and anti self-dual forms.
In higher dimensions the set ${\cal S}_{2n}$ is an $n^2-n+1$ dimensional
submanifold [12].

We recall the inequalities
$$(\omega,\eta)^2\le(\omega,\omega)(\eta,\eta),\quad\quad
2(\omega,\eta)\le (\omega,\omega)+(\eta,\eta).\eqno(2.3)$$
From the (2.1) and (2.3) we obtain the following Lemma which
generalizes the results of [8] to arbitrary dimensions.

\proclaim Lemma 2.2. Let $\omega$ be a 2-form in $2n$ dimensions. Then
$$(n-1)
(\omega,\omega)^2-{n\over 2}(\omega^2,\omega^2)\ge0\eqno(2.4)$$
$$(\omega^{n/2},\omega^{n/2})\ge *\omega^n,\eqno(2.5)$$
and equality holds if and only if all eigenvalues of $\omega$ are equal.

\vskip .2cm
\noindent
{\it Proof.} To obtain  Eq.(2.4) we use,
$$\eqalignno{
q_1^2\ge &q_2\cr
{1\over n^2}s_2^2\ge&{2\over n(n-1)}s_4\cr
{1\over n^2}(\omega,\omega)^2\ge&{2\over n(n-1)}{1\over
4}(\omega^2,\omega^2)\cr}$$
which gives the desired result. Similarly using
$$\eqalignno{
q_{n/2}^2\ge & q_n\cr
\left({(n/2)!(n/2)!\over n!}s_{n/2}\right)^2\ge&s_{2n}\cr
\left({1\over n!}(\omega^{n/2},\omega^{n/2})\right)^2\ge&
{1\over (n!)^2}\mid *\omega^n\mid^2\cr}$$
we obtain (Eq.(2.5).\quad \quad e.o.p.

From Lemma 2.2,  we immediately have

\proclaim Corollary 2.3. The 2-form $\omega$ is strongly self-dual iff
$\omega^{n/2}$ is self-dual in the Hodge sense.

We will now show that this condition is
also equivalent to the self-duality definition used by Trautman [3].

\proclaim Proposition 2.4. Let $\omega $ be a 2-form in $2n$ dimensions. Then
$$\omega^{n-1}=k*\omega\eqno(2.6)$$ where $k$ is a constant, if and only if
$\omega$ is strongly self-dual and
$k={n!\over n^{n/2}}
(\omega,\omega)^{{n\over 2}-1}$.

\noindent
{\it Proof.}
Multiplying (2.6) with $\omega$ and taking Hodge duals,
we obtain,
$*\omega^n=k(\omega,\omega)$. Since $(\omega,\omega)=s_2=nq_1$ and
$\mid*\omega^n\mid=n!s_{2n}^{1/2}=n!q_n$, we obtain $k=(n-1)!q_n^{1/2}/q_1$.
Then
taking inner products
of both sides of (2.6) with themselves, we obtain
$(\omega^{n-1},\omega^{n-1})=k^2(*\omega,*\omega)=k^2(\omega,\omega)$.
Substituting the vakue of $k$ obtained above,  and using
$(\omega^{n-1},\omega^{n-1})=\big((n-1)!\big)^2nq_{n-1}$, we obtain
$q_n= q_{n-1}q_1$. But since $q_1\ge q_n^{1/n}$, we have $q_n\ge
q_{n-1}q_n^{1/n}$, which leads to $q_n^{n-1}\ge q_{n-1}^n$. This is just
the reverse of the inequlity in  (2.2), hence equality must hold,
and all eigenvalues of $\omega$ are equal in absolute value.
Thus  $\omega$ is strongly self-dual and
it can also  be seen that
$k={n!\over n^{n/2}}
(\omega,\omega)^{{n\over 2}-1}$. \quad e.o.p.
\vskip .3cm

Applying Lemma 2.2 to the 2-forms $\omega\pm\eta$ we obtain

\proclaim Corollary 2.5. Let $\omega $ and $\eta$ be 2-forms in $2n$
dimensions. Then
$$\eqalignno{
4(\omega\eta,\omega\eta)\le &4\t{2(n-1)\over n}(\omega,\eta)^2
+\left[\t{2(n-1)\over n)}(\omega,\omega)^2-(\omega^2,\omega^2)\right]
+\left[\t{2(n-1)\over n}(\eta  ,\eta  )^2-(\eta  ^2,\eta  ^2)\right]\cr
&+2\left[\t{2(n-1)\over
n}(\omega,\omega)(\eta,\eta)-(\omega^2,\eta^2)\right]}$$

The inequality is saturated when $\omega=\eta$ and it is strongly self-dual.
This corollary will be used in the next section.

\vskip 1cm

\noindent
{\bf 3. An upper bound bound for $\sigma_4$.}

Consider an $SO(N)$ valued curvature 2-form $F=(F_{ab})$. Then
$$\eqalignno{
\langle F, F\rangle=&2\sum(F_{ab},F_{ab}),&(3.1a)\cr
\sigma_2=&\sum F_{ab}^2&(3.1b),\cr
\sigma_4=&\sum\left[F_{ab}F_{cd}-F_{ac}F_{bd}+F_{ad}F_{bc}\right]^2,&(3.1c)\cr}
$$
where the the summation goes over the indices such that $b>a$, $d>c$, and in
the last summation the the pairs $\{a,b\}$ and $\{c,d\}$ are distinct with
$c>a$. We recall the integrals  $\int (\sigma_2)^k$ (over a $4k$ manifold)
and  $\int (\sigma_4)^k$ (over an $8k$ manifold) are
independent of the bundle trivilization {\it and} of the representatives of the
cohomology classes, hence they are invariants
of the bundle, while, $\int\langle F,F\rangle^k$ and $\int(\sigma,\sigma)^k$
are
only independent of the bundle trivilizations. Thus our aim is to obtain upper
bounds
for the bundle invariant quantities in terms of latter quantities, which can be
taught of physical fields, and to characterize the conditions under which these
inequalities are saturated.

In [8] we have obtained an upper bound for $\int \sigma_2^2$  for an $SO(N)$
bundle over an eight dimensional manifold. This result can be generalized
immediately to $SO(N)$ bundles over $2n$ dimensional manifolds.

\proclaim Proposition 3.1. Let $F$ be the curvature 2-form of an
$SO(N)$ bundle over a $2n$ dimensional
manifold,
Then
  $$\langle F,F\rangle^2\ge \t {1\over 4}{2(n-1)\over n}
(\sigma_2,\sigma_2)\ge\t {1\over 4}{2(n-1)\over n}   *\sigma_2^2.\eqno(3.2)$$

The conditions under which the inequalities are saturated depend both on the
relations among the $F_{ab}$'s,  and on the properties of each $F_{ab}$, i.e.,
inequalities of type (2.3) and (2.4) should be both saturated. For example,
(3.2) is saturated when $F=\omega F_o$, where $F_o$ is a constant matrix, and
$\omega$ is a strongly self-dual 2-form. By raising both sides of the
inequality in (3.2) to power $k$ and integrating, we obtain upper
bounds for the integral of $p_1^k$ on a $4k$ dimensional manifold.

We can give estimates on $\sigma_4$ for an $SO(N)$ bundle using Corollary 2.5.
To illustrate the proof, we will first work for an $SO(4)$ bundle, then
generalize the result to arbitrary dimensions.

\proclaim Proposition 3.2. Let $F$ be the curvature 2-form of an $SO(4)$
bundle over
an eight dimensional manifold. Then
$$\eqalignno{\mid *\sigma_4\mid
\le &\t {9\over 2}\Phi +\t{9\over 32}\langle F,F\rangle^2-{3\over 4}
(\sigma_2,\sigma_2)&(3.3)\cr}$$
where $\Phi=(F_{12},F_{34})^2+(F_{13},F_{24})^2+(F_{14},F_{23})^2$.

\noindent
{\it Proof.} We explicitly write
$$\langle F,F\rangle=2\left[
(F_{12},F_{12})+(F_{13},F_{13})+(F_{14},F_{14})+(F_{23},F_{23})+(F_{24},F_{24})
+ (F_{34},F_{34})\right],$$
$$\sigma_2=
F_{12}^2+F_{13}^2+F_{14}^2+F_{23}^2+F_{24}^2+F_{34}^2,$$
$$\sigma_4=\left[F_{12}F_{34}-F_{13}F_{24}+F_{14}F_{23}\right]^2.$$
Then using the inequality  $*\eta^2<(\eta,\eta)$ for 4-forms, we
have
 $$\mid *\sigma_4\mid\le
 (F_{12}F_{34}-F_{13}F_{24}+F_{14}F_{23},
 F_{12}F_{34}-F_{13}F_{24}+F_{14}F_{23}),$$
which gives using (2.3b)
 $$\mid *\sigma_4\mid\le
3\left[(F_{12}F_{34},F_{12}F_{34})+
            (F_{13}F_{24},F_{23}F_{24})+
            (F_{14}F_{23},F_{14}F_{23})\right]$$
Now we can use Corollary 2.5
to write the right hand side in terms of
$\langle F,F\rangle$ and $(\sigma_2,\sigma_2)$. However
we need to add
terms such as
$3(F_{12},F_{12})(F_{13},F_{13})-2(F_{12}^2,F_{13}^2)$
(which are all positive quantities) in order to obtain $\langle F,F\rangle$
and $(\sigma_2,\sigma_2)$ on the right hand side. Collecting terms
we obtain (3.3).
 \quad e.o.p.

 \vskip .3cm
The inequality is saturated when $F_{12}=F_{34}$, $F_{13}=-F_{24}$,
$F_{14}=F_{23}$, and each of these 2-forms are strongly self-dual.
The quantity $\Phi$ is not invariant under bundle automorphisms in general,
however it can be seen that the condition $\Phi=0$ is invariant if
$$(F_{ij},F_{kl})=\delta_{ik}\delta_{jl}\eqno (3.4).$$
It is possible to obtain
a bound for
$\Phi$ in terms of $\langle F,F\rangle$, but then the inequality (3.3) is not
saturated.

The Proposition can be generalized to arbitrary dimensions by counting
arguments.

\proclaim Corollary 3.3. Let $F$ be the curvature 2-form of an $SO(N)$ bundle
over a $2n$ dimensional  manifold. Then
$$\mid *\sigma_4\mid\le\t{3\over 4}\left[ 6\Phi
+\t{3\over 8}{N-2\choose 2}\langle F,F\rangle^2-
{N-2\choose 2}(\sigma_2,\sigma_2) \right ]\eqno(3.5)$$
where $\Phi=\sum(F_{ij},F_{kl})$, where the summation extends over
distinct pairs of indices $\{i,j\}$ and $\{k,l\}$.

By raising both sides of (3.4) to power $k$ and integrating we obtain upper
bounds for $\int (p_2)^k$ over an $8k$ dimensional manifold.
The right hand side of Eq.(3.5) is gauge invariant only if
$\Phi=0$ and an analogue of Eq.(3.4) holds.

\vskip 1cm
\noindent
{\bf 4. Conclusion.}

The equivalence of strong self-duality and various other
definitions of self-duality in dimensions higher than four, and the saturation
of certain topological lower
bounds suggests that strongly self-dual forms share many important properties
of self-dual 2-forms in four dimensions.  There are however certain features
that are not  matched completely.

First of all
the strongly self-dual 2-forms do not form a linear space. This situation
can be remedied by restricting  the notion of self-duality to linear
submanifolds [12]. However the maximal dimension of such linear submanifolds in
dimensions $2n$ is equal to the number of linearly independent vector fields on
$S^{2n-1}$. In eight dimensions those linear submanifolds include a
representation of octonions in $R^8$ and the self-duality equations of Corrigan
et al [2]. However in dimensions for example $2(2a+1)$ these linear subspaces
are one dimensional, and the self-duality equations obtained in [2] cannot be
represented in this framework. This suggests that the notion of strong
self-duality may be weekened to comprise a richer structure of self-duality
equations, at the price of relaxing the saturation of topological lower bounds.

 The second problem is that, the  self-duality  equations obtained by selecting
linear subspaces of strongly self-dual forms  are overdetermined. Thus
the existence
of solutions are not guaranteed in general (one solution is given in [6]). We
remark however that the in eight dimensions, the orthogonal complements of the
linear subspaces of strongly self-dual forms lead to elliptic equations, and
this suggest that in arbitrary dimensions, these complememts may be more
interesting from the point of view of resulting system of partial differential
equations.

\vskip 1cm
\noindent
{\bf Acknowledgements.} This work is partially supported by the Turkish
Scientific and Technological Research Council, TUBITAK.
\vskip 1cm

\baselineskip 16pt

\noindent
 {\bf References.}
\vskip .3cm
\item{[1]}F.W.Warner, {\it Foundations of Differentiable Manifolds and Lie
Groups}, Springer, New York, (1983).

\item{[2]}E.Corrigan, C.Devchand, D.B.Fairlie and J.Nuyts,
{\it
Nuclear Physics}, {\bf B214}, 452-464, (1983).

\item{[3]}A.Trautman, {\it International Journal of Theoretical
Physics}, {\bf 16}, 561-565, (1977).

\item{[4]}D.H. Tchrakian, {\it Journal of Mathematical Physics}, {\bf
21}, 166, (1980).

\item{[5]}C. Saclioglu, {\it Nucleaar Physics}, {\bf B277}, 487,
(1986).

\item{[6]}B.Grossman, T.W.Kephart, J.D.Stasheff,
{\it Commun. Math.
Phys.}, {\bf 96}, 431-437, (1984).

\item {[7]} R.S. Ward, {\it Nucl. Phys.}, {\bf B236}, 381, (1984).

\item{[8]} A.H. Bilge, T. Dereli and S. Kocak, ``Self-dual Yang-Mills fields
in eight dimensions", {\it Lett. Math. Phys.}, to appear.

\item{[9]}J.M.Milnor and J.D. Stasheff, {\it Characterisitic Classes}, Annals
of
Mathematics Studies, Princton University Press, Princton, New Jersey, (1974).

\item{[10]}F.R.Gantmacher, {\it The Theory of Matrices}, Chelsea Publ.
Co., New York, (1960).

\item{[11]}M. Marcus and H.Minc, {\it A Survey of Matrix Theory and
Matrix Inequalities}, Allyn and Bacon Inc., Boston, (1964).

\item{12]} A.H. Bilge, T. Dereli and S. Kocak, ``An explicit
construction of self-dual 2-forms in eight dimensions", hep-th
9509041.

\end